\documentclass[12pt]{amsart}
\usepackage{color}
\usepackage{amssymb}
\usepackage{graphicx}
\usepackage{ mathrsfs }

\newtheorem{Theorem}{Theorem}

\begin{document}

\title[THE RING OF CONDITIONS FOR  $(\Bbb C^*)^n$]{NEWTON POLYHEDRA, TROPICAL GEOMETRY AND THE RING OF CONDITIONS FOR  $(\Bbb C^*)^n$ }

\author{Boris Kazarnovskii}
\address{Institute for Information Transmission Problems (Kharkevich Institute), Russian Academy
of Sciences, Bolshoy Karetny per. 19, build. 1, Moscow 127051, Russia}
\email{kazbori@gmail.com}

\author{Askold Khovanskii}
\address{Department of Mathematics, University of Toronto, Toronto,
Canada; Moscow Independent University, Moscow, Russia.}
\email{askold@math.toronto.ca}

\begin{abstract}
The  ring of conditions defined by C. De Concini and C. Procesi  is an  intersection theory  for
algebraic cycles in a spherical homogeneous space. In the paper we consider the ring of conditions
for the group $(\Bbb C^*)^n$. Up to a big extend this ring  can be reduced to the cohomology rings
of smooth projective toric varieties. This ring also can be described using tropical geometry.
We recall these known results and   provide a new description  of this ring in terms of convex
integral polyhedra.

\end{abstract}
\thanks{ The first author was partially supported by the Russian Foundation of Sciences, Project No. 14-50-00150. The second author was partially supported by the Canadian Grant No. 156833-17.}

\date{\today}
\maketitle

\section{Introduction}

The  ring of conditions (see~\cite{Conc_Pro}) is an  intersection theory  for algebraic cycles in a spherical homogeneous space  with coefficients in a commutative ring $\Lambda$.  In the paper we consider the ring of conditions  $\mathcal{R}_n (\Lambda)$ for the group $(\Bbb C^*)^n$ (which is a spherical homogeneous space with respect to the natural action of the group on itself)  with coefficients in  $\Lambda=\Bbb Z, \Bbb R, \Bbb C$. Up to a big extend the ring $\mathcal{R}_n (\Lambda)$  can be reduced to the cohomology rings of smooth projective toric varieties \cite{Conc_Pro}. Tropical geometry,  relates  algebraic geometry and piecewise linear geometry (see~~\cite{Ful_Stu} -- \cite{Ma_Stu}). In particular, it studies tropicalization of subvarieties in $(\Bbb C^*)^n$ and their intersection theory (see~\cite{Ful_Stu}). We   remind these results and  provide a new description  of $\mathcal {R}_n(\Lambda)$ in terms of convex integral polyhedra. Thus we  show that in fact the ring $\mathcal {R}_n(\Lambda)$ belongs to  Newton polyhedra theory.

\section{The ring  of conditions  $\mathcal{R}_n (\Lambda)$ for the group $(\Bbb C^*)^n$}

Two $k$-dimensional cycles $X_1,X_2\subset (\Bbb C^*)^n$  are {\it equivalent} $X_1\sim X_2$ if for any $(n-k)$-dimensional  cycle  $Y\subset (\Bbb C^*)^n$ and for
almost all  $g\in (\Bbb C^*)^n$ we have $ <X_1, gY>= <X_2, gY>$ where $<A,B>$ is the intersection index of the cycles $A$ and $B$. If $X_1\sim X_2$ and $Y_1\sim Y_2$ then for almost all $g_1,g_2\in (\Bbb C^*)^n$ we have $X_1\cap g_1Y_1 \sim X_2\cap g_2Y_2$. {\it The product} $X\ast Y$ of  equivalence classes $X$ and $Y$ is the equivalence classes of the intersection $X_1\cap g_1Y_1$ where $X_1$ and $Y_1$ are representatives of $X$ and $Y$ and $g_1$ is a generic element in $(\Bbb C^*)^n$. {\it The ring of conditions} $\mathcal{R}_n (\Lambda)$  is
the ring of  the equivalence classes  with the  multiplication $\ast $ and with the tautological addition.

\section {Bergman cone}

A vector $k=(k_1,\dots,k_n)\in \Bbb Z^n$ is  {\it essential} for a  variety $X\subset (\Bbb C^*)^n$ if there is a germ of meromorphic map $f:(\Bbb C,0)\rightarrow X\subset (\Bbb C^*)^n$ where  $f(t)=at^k+\dots$,  $a=(a_1,\dots,a_n)\in (\Bbb C^*)^n$ and points $\dots$ denote terms $a_mt^m$ of higher order (i.e. $m=(m_1,\dots,m_n)$ where $m_i\geq k_i$ for $1\leq i\leq n$ and $m\neq k$).

The {\it Bergman cone} $B(X)\subset \Bbb R^n$ of $X$ is the closure of the set of vectors $\lambda k\in \Bbb R^n$  where $k$ is  essential vector for $X$ and $\lambda\geq 0$.
\begin{Theorem} If each irreducible component of $X$ has complex dimension $m$ then $B(X)$ is a finite union of convex rational cones $|\sigma_i|\subset \Bbb R^n$ with $\dim _{\Bbb R}|\sigma_i|=m$. Moreover $B(X)$ is the support of a fan (defined up to subdivision of $B(X)$) of some toric variety.
\end{Theorem}
A first version of theorem 1 appeared in \cite{Ber}. Now its different versions can be found in many works dedicated to tropical geometry.

\section {Good compactification}

 Toric variety $M\supset (\Bbb C^*)^n $ is a {\it good compactification} for a  subvariety $X\subset (\Bbb C^*)^n$ with $\dim X=k$ if its closure $\overline X$ in $M$  is complete and does not intersect orbits in $M$  whose codimension is bigger than $k$.

\begin{Theorem} 1) For any finite set $\mathcal{S}$ of algebraic subvarieties in $(\Bbb C^*)^n$ there is a toric variety $M\supset (\Bbb C^*)^n$ which provides a good compactification for each subvariety from $\mathcal{S}$.

2) Toric variety $M$ is a good compactification of $X\subset (\Bbb C^*)^n$ if and only if the support of its fan contains the Bergman cone $B(X)$.
\end{Theorem}

The part 1) of  theorem 2 was discovered (in a stronger form) in  \cite{Conc_Pro}. It is essential for the theory of rings of conditions. Now  different versions of theorem 2 can be found in many works dedicated to tropical geometry (for example, see~\cite{Tev}). Let $\mathcal{S}_r$ be a set of all subvarieties in $(\Bbb C^*)^n$ such that any
$X$ from $\mathcal{S}_r$ can be defined by a system of Laurent polynomials whose Newton polyhedra belong to a ball of radius $r$.The following  more precise version of the part 1) of theorem 2 easily follows from  \cite{Ka_Kh I}.
\begin{Theorem} There is a Newton polyhedron $\Delta_r$ such that the  projective toric variety $M_{\Delta_r}$ corresponding to $\Delta_r$ is smooth and it provides a good compactification for any $X\in \mathcal{S}_r$. Bergman set $B(X)$ of any $X\in \mathcal{S}_r$ is a subfan of the dual fan $\Delta_r^\bot$ to the polyhedron $\Delta_r$.
\end{Theorem}

\section {The ring  $\mathcal{R}_n (\Lambda)$ and cohomology rings of toric varieties}

For a complete smooth toric variety $M\supset (\Bbb C^*)^n$ and for any $k$-dimensional cycle $X=\sum k_iX_i$ one can defined the cycle $\overline X$ in $M$ as $\sum k_i\overline X_i$ where $\overline X_i$ is the closure  in $M$ of  $X_i\subset (\Bbb C^*)^n$. The cycle $\overline X$  defines an element $\rho (\overline X)$ in $H^{2(n-k)}(M^n, \Lambda)$ whose value on the closure $\overline O_i$ of  an $(n-k)$-dimensional orbit $O_i$ in $M$ is equal to the intersection index $<\overline X,\overline O_i>$. A  compactification $M\supset (\Bbb C^*)^n $ is {\it good} for a cycle $X=\sum k_iX^i$ in $(\Bbb C^*)^n$ if it is good compactification for each $X_i$.
\begin{Theorem}[see \cite{Conc_Pro}] If a smooth toric compactification $M$ is good for  cycles $X, Y$ and $Z$ where $Z=X\ast Y$, then the product $\rho(X)\rho(Y)$ in the cohomology ring $H^*(M,\Lambda)$ of the elements $\rho(X)$ and $\rho(Y)$  is equal to $\rho (Z)$.
\end{Theorem}

\section {The ring of balanced $\Lambda$-enriched fans}

{\bf 6.1.} An {\it $\Lambda$-enriched $k$-fan}  is a fan $\mathcal F\subset \Bbb R^n$ of an $n$-dimensional toric variety equipped  with a {\it weight function}  $c:\mathcal{F}_k \rightarrow \Lambda$   defined on the set $\mathcal{F}_k$  of all $k$-dimensional cones from $\mathcal{F}$. The {\it support} $|\mathcal F|$ of  $\mathcal F$ is the union of  all  cones $|\sigma_i|\subset \Bbb R^n$ where   $\sigma_i\in\mathcal {F}_k$ and $c(\sigma_i)\neq 0$. Two enriched $k$-fans $\mathcal{F}_1$ and $\mathcal{F}_2$ are {\it equivalent} if:
1) $|\mathcal{F}_1|=|\mathcal{F}_2|$; 2) the weight functions $c_1$ and $c_2$ induce the same weight function on every common subdivision of the fans $\mathcal{F}_1$ and $\mathcal{F}_2$.

{\bf 6.2.}
Let $\mathcal{F}$ be an enriched $k$-fan. For a cone $\sigma_i\in \mathcal{F}_k$ let $L_i^\bot\subset (\Bbb R^n)^*$ be the space  dual to  the span $L_i$ of $|\sigma_i|\subset \Bbb R^n$. Let $O$ be an orientation  of $|\sigma_i|$.  Denote by $e_i^\bot(O)\in \Lambda^{n-k}L_i^\bot$  the $(n-k)$-vector,  such that: 1) the integral volume of $|e_i^\bot(O)|$  in $L_i^\bot$ is equal to one; 2) the orientation of $e_i^\bot(O)$ is induced from the orientation $O$ of $|\sigma_i|$ and from the standard orientation of  $\Bbb R^n$. An enriched  $k$-fan $\mathcal{F}$ satisfies {\it the balance condition} if for any orientation of a $(k-1)$-dimensional cone $|\rho|$ where $\rho\in F_{k-1}$, the  relation
\begin{equation}
\label{balance}
\sum e_i^\bot (O(\rho))c(\sigma_i)=0
\end{equation}
holds, where $c$ is the weight function and summation  is taken over all $\sigma_i\in \mathcal{F}_k$ such that $|\rho|\subset \partial |\sigma_i|$ and $O(\rho)$ is such  orientation of $|\sigma_i|$ that the orientation of $\partial |\sigma_i|$ agrees with the orientation of $|\rho|$.

{\bf 6.3.} Let $\mathcal{F}_1$ and  $\mathcal{F}_2$ be balanced $k$- and  $(n-k)$-fans.
Cones $\sigma^1_i\in \mathcal{F}_1$, $\sigma^2_j\in \mathcal{F}_2$ with $\dim \sigma_i^1=k$, $\dim \sigma_j^2=n-k$ are {\it $a$-admissible} for a vector $a\in \Bbb R^n$ if $|\sigma^1_i|\cap (|\sigma^2_j|+a)\neq \emptyset$. Let $C_{i,j}$ be the index of $\Lambda_i\bigoplus\Lambda_j$ in $\Bbb Z^n$ where $\Lambda_i=L_i^1\cap\Bbb Z^n$, $\Lambda_j=L_j^2\cap\Bbb Z^n$ and $L_i^1$, $L_j^2$ are linear spaces spanned by $|\sigma^1_i|$, $|\sigma^2_j|$. The {\it intersection number} $c(0)$ of $\mathcal{F}_1$and $\mathcal{F}_2$ is
\begin{equation}
\label{sum}
\sum C_{i,j}c_1(\sigma^1_i)c_2(\sigma^2_j),
\end{equation}
where the sum is taken over all $a$-admissible couples $\sigma^1_i,\sigma^2_j$ for any generic vector $a\in \Bbb R^n$ (one can show that if  the fans $\mathcal{F}_1$,  $\mathcal{F}_2$ satisfy the balance condition (\ref{balance}) then the sum (\ref{sum}) is independent of the choice of generic vector $a$).
The {\it product} $\mathcal{F}=\mathcal{F}_1\times \mathcal{F}_2$ is a $0$-fan $\mathcal{F}=\{0\}$ with the weight $c(0)$ equal to the intersection number.

{\bf 6.4.}
Consider a $k$-fan $\mathcal{F}_1$ and a $m$-fan $\mathcal{F}_2$  from the set $T\mathcal{R}_n (\Lambda)$ of all balanced $\Lambda$-enriched fans. Let $d$ be $n-(k+m)$.  If $d<0$ then $\mathcal{F}_1\times\mathcal{F}_2 =0$. If $d=0$ the fan $\mathcal{F}_1\times\mathcal{F}_2$ is  defined above. Let us define the $d$-fan $\mathcal{F}=\mathcal{F}_1\times\mathcal{F}_2$ for $d>0$. Assume that $\mathcal{F}_1$ and and  $\mathcal{F}_2$ are subfans of a complete fan $\mathcal{G}$. Then $\mathcal{F}=\mathcal{F}_1\times\mathcal{F}_2$ also is a subfan of $\mathcal{G}$. The weight $c(\delta)$ of a cone $\delta$ from $\mathcal{G}$ with $\dim \delta=d$   is defined as follows. Let $L$ be a space spanned by the cone $|\delta|$ and let $(\mathcal{F}_1)_\delta$ and $(\mathcal{F}_2)_\delta$ be the enriched subfans of $\mathcal{F}_1$ and of $\mathcal{F}_2$ consisting of all cones from these fans containing the cone $\delta$. The  weight $c(\delta)$ of the cone $\delta$ in $\mathcal{F}=\mathcal{F}_1\times\mathcal{F}_2$ is equal to the intersection number of the images under the factorization of $(\mathcal{F}_1)_\delta$ and $(\mathcal{F}_2)_\delta$ in the factor space $\Bbb R^n/L$ equipped with the factor lattice $\Bbb Z^n/ L\cap \Bbb Z^n$.

\section{Tropicalization of the ring  $\mathcal{R}_n (\Lambda)$}

Let $\Delta^\bot$ be a fan of a smooth complete projective toric variety $M^n_\Delta$. Let $T\mathcal{R}_n(\Lambda,\Delta)$ be a ring of balanced $\Lambda$-enriched fans equal to $\Lambda$-linear combination of cones from the fan $\Delta^\bot$. The following theorems 5,6 are proved in \cite{Ful_Stu}.

\begin{Theorem}[see \cite{Ful_Stu}] The ring $T\mathcal{R}_n(\Lambda,\Delta)$ is isomorphic to the intersection ring $H_{*}(M_\Delta,\Lambda)$. The component of $T\mathcal{R}_n(\Lambda,\Delta)$ consisting of $k$-fans under this isomorphism corresponds to the component $H_{2k}(M_\Delta,\Lambda)$.
\end{Theorem}

\begin{Theorem}[see \cite{Ful_Stu}]\label{theorem 5} The ring of conditions $\mathcal{R}_n(\Lambda)$ is isomorphic to the tropical ring $T\mathcal{R}_n(\Lambda)$  of all balanced $\Lambda$-enriched fans.(The rings $\mathcal{R}_n(\Bbb Z)$, $\mathcal{R}_n(\Bbb C)$ have  similar descriptions).
\end{Theorem}

\section{The graded ring associated to a homogeneous polynomial}

To a homogeneous polynomial $P$ on a   $\Bbb R$-linear space $\mathcal {L}$, $\dim \mathcal {L}<\infty$, one can associate the graded commutative ring $A (\mathcal {L},P)$ (one can produce a similar constructions for  homogeneous  polynomials on  infinite dimensional spaces over any field and for functions  analogues to homogeneous  polynomials  on free abelian groups). Let $D(\mathcal{L})$ be the ring of $\Bbb R$-linear differential operators on $\mathcal{L}$ with constant  coefficients. It is generated by Lie derivatives $L_v$ along  constant vector fields $v(x) \equiv v\in \mathcal {L}$ and by multiplications on real constants. Let $I_P\subset D(\mathcal{L})$ be  a set defined by the following condition: $L\in I_P\Leftrightarrow L (P) \equiv 0$. It is easy to see that $I_P$ is a homogeneous ideal. By definition the {\it ring associated to $P$ } is the  factor ring $A (\mathcal {L},P)=D(\mathcal {L}) /I_P$. One can to see that: (1) $A(\mathcal{L},P)$ is a graded  ring with homogeneous  components $A_0,\dots, A_n$ where $n=\deg P$; (2)  $A_0= \Bbb R$;
(3) there is a non-degenerate pairing between $A_k$ and $A_{n-k}$ with values in $A_0$, thus $ A_k= A_{(n-k)}^*$ and $A_n\sim \Bbb R$.

\section{Khovanskii--Pukhlikov ring and the ring  $\mathcal{R}_n (\Lambda)$}

Let $L_\Delta$ be a  space of formal differences of convex polyhedra whose support functions are linear on each cone from the fan of a smooth projective toric variety $M_\Delta$. Let $n!V$ be the degree $n$ homogeneous polynomial on $L_\Delta$  whose  value  on  $\tilde \Delta\in  {L}_\Delta$ is equal to the volume of $\tilde \Delta$ multiplied by $n!$. The ring $A(L_\Delta, n!V)$ is called Khovanskii--Pukhlikov ring \cite{Pu_Kh}.

\begin{Theorem} The intersection ring  $H_*(M_\Delta,\Bbb R)$  is isomorphic (up to a change of the grading) to the Khovanskii--Pukhlikov ring  $A(L_\Delta, n!V)$.
\end{Theorem}

Let $\mathcal{L}_n$, $\dim \mathcal{L}_n =\infty$, be the space of formal differences of convex polyhedra $\Delta$ with rational dual fans $\Delta^\bot$. Let  $n!V$ be the degree $n$ homogeneous polynomial on $\mathcal{L}_n$  whose value  on  $\Delta\in \mathcal {L}_n$ is equal to  the volume of $\Delta$ multiplied by n!.

\begin{Theorem}\label{theorem 7} The ring $\mathcal{R}_n(\Bbb R)$ is isomorphic to the ring  $A(\mathcal{L}_n, n!V)$. (The rings $\mathcal{R}_n(\Bbb Z)$, $\mathcal{R}_n(\Bbb C)$ have  similar descriptions).
\end{Theorem}

\section{The BKK theorem and the ring $\mathcal{R}_n (\Lambda)$}

Let $\{\Gamma_i\}$ be a set of $n$ hypersurfaces in $(\Bbb C^*)^n$  defined by $P_i=0$ where $P_i$ are Laurent polynomials with Newton polyhedra $\Delta_i$. Bernstein-Koushnirenko-Khovanskii theorem (BKK theorem \cite{Bern}-\cite{Kh} ) can be stated in the following two ways:
\begin{Theorem} The intersection number of the hypersurfaces $\Gamma_i$ in the ring of conditions is equal to the mixed volume of $\Delta_1,\dots,\Delta_n$ multiplied by $n!$.
\end{Theorem}

Let $\mathcal{F}_i$ be $\Bbb R$-enriched $(n-1)$-fan dual to  $\Delta_i$ whose weight function  at a cone $\sigma$ dual to a side $\sigma^\bot$ of $\Delta_i$ is equal to the integral length of the $\sigma^\bot$.

\begin{Theorem} The intersection number of the hypersurfaces $\Gamma_i$ in the ring of conditions is equal to the intersection number of the $\Bbb R$-enriched fans $\mathcal{F}_i$  in the ring~$T\mathcal{R}_n$.
\end{Theorem}

Thus theorems \ref{theorem 7} and \ref{theorem 5} could be considered as  generalizations of the BKK theorem. Such generalizations are possible because of the following reason: the cohomology ring of a smooth toric variety is generated by elements of degree two.


\begin{thebibliography}{99} 

\bibitem{Conc_Pro}
C.~De Concini and C.~Procesi, \textit{Complete symmetric varieties II Intersection
theory.}, Adv. Stud. Pure Math., \textbf{6} (1985), 481--512.
\bibitem{Ful_Stu}
W.~Fulton, B.~Sturmfels, \textit{Intersection theory on toric varieties}, Topology, \textbf{6}, No.~2, (1997), 335--353.

\bibitem{Kaz}
B.~Kazarnovskii,\textit {Truncations of systems of equations, ideals and
varieties}, Izv. Math.\textbf{63}, No.~3, (1999),  535--547.


\bibitem{Iten_Mikh_Shu}
I.~Itenberg, G.~Mikhalkin, E.~Shustin, \textit{Tropical algebraic geometry}, (2nd ed.), Birkhäuser, Basel(2009).

\bibitem{Ma_Stu}
D.~Maclagan, B.~Sturmfels, \textit{Introduction to Tropical Geometry},
(Graduate Studies in Mathematics), (2015).

\bibitem{Ber}
G.~Bergman,\textit {The logarithmic limit-set of an algebraic variety}, Trans. Amer. Math. Soc. \textbf {157} (1971), 459--469.
\bibitem{Tev}
J.~Tevelev, \textit {Compactifications of  subvarieties of tori}
American Journal of Mathematics
\textbf {129}, No.~4 (2007), 1087--1104.

\bibitem{Ka_Kh I}
 B.~Kazarnovskii and A.~Khovanskii, \textit {Tropical Notherian property and Gr¨obner bases,} St. Petersburg Math. J.  \textbf {26} (2015), No. 5,  797--811.

\bibitem{Pu_Kh}
\textit {A.~Pukhlikov and A.~Khovanskii,} \textit{A Riemann–-Roch theorem for integrals and sums of quasipolynomials over virtual polytopes,}  St. Petersburg Math. J.  \textbf{4} (1993),No.~4,  789--812.

\bibitem{Bern}
D.~Bernstein, \textit{The number of roots of a system of equations,} Funkts. Anal. Prilozhen. \textbf {9}, (1975), No.~3, 1--4.



\bibitem{Kou}
A.~Kouchnirenko, \textit {Polyedrea de Newton et nombres de Milnor,} Inv. Math. \textbf {32}, 1--31.

\bibitem {Kh}
A.~Khovanskii,\textit { Newton polyhedron, Hilbert polynomial and sums of finite sets,}  Funkts. Anal. Prilozhen.\textbf {26}, (1993), No. 4, 276--281.

\end{thebibliography}
\end{document}